\documentclass[12pt,a4paper]{amsart}

\usepackage {amsfonts}
\usepackage[english]{babel}
\def\b{\beta}
\def\g{\gamma}
\def\G{\Gamma}
\def\d{\delta}
\def\D{\Delta}
\def\a{\alpha}
\def\p{\varphi}
\def\e{\varepsilon}
\def\l{\lambda}
\def\L{\Lambda}

\def\la{\langle}
\def\ra{\rangle}
\def\R{{\mathbb R}}
\def\C{{\mathbb C}}
\def\N{{\mathbb N}}
\def\Z{{\mathbb Z}}

\def\p{\varphi}
\def\bs{~\hfill\rule{7pt}{7pt}}
\def\U{\Upsilon}

\DeclareMathOperator{\supp}{supp}

\newtheorem{Th}{Theorem}

\newtheorem*{Cor}{Corollary}

\newtheorem{Pro}{Proposition}

\begin{document}

\title{Local Wiener's Theorem and Coherent Sets of Frequencies}

\author{S.Yu.Favorov}

\address{Serhii Favorov,
\newline\hphantom{iii}  Karazin's Kharkiv National University
\newline\hphantom{iii} Svobody sq., 4,
\newline\hphantom{iii} 61022, Kharkiv, Ukraine}
\email{sfavorov@gmail.com}

\maketitle {\small
\begin{quote}
\noindent{\bf Abstract.}
Using a local analog of the Wiener-Levi theorem, we investigate the class of  measures on Euclidean space with  discrete support and spectrum. Also, we find a new sufficient conditions for a discrete set in Euclidean space to be a coherent set of frequencies.
\medskip


\end{quote}
}

\medskip
{\bf 1. The Wiener--Levi Theorem}. The following theorem is known as the Wiener--Levi Theorem (see, for example, \cite{Z}, Ch.VI):
\begin{Th}
Let
$$
F(t)=\sum_{n\in\Z}c_ne^{2\pi int}
$$
 be an absolutely convergent Fourier series,  and $h(z)$  be a holomorphic function on a neighborhood of the closure of the set $\{F(t):\,t\in [0,1]\}$. Then the function $h(F(t))$ admits  an absolutely convergent Fourier series expansion as well.
 \end{Th}
 For $h(z)=1/z$ the theorem is well-known as Wiener's Theorem.
 \medskip

 Denote by $W$ the algebra of absolutely convergent series
 $$
 f(x)=\sum_n c_ne^{2\pi i\la x,\l_n\ra},\quad \l_n\in\R^d,\  x\in\R^d,
 $$
 with the norm $\|f\|_W=\sum_n|c_n|$. In \cite{F0} we proved
 \begin{Th}\label{T}
 Let  $K\subset\C$ be an arbitrary compact set,  $h(z)$ be a holomorphic function on a neighborhood of $K$, and $f\in W$. Then there is a function $g\in W$ such that if $f(x)\in K$ then $h(f(x))=g(x)$.
 \end{Th}
 For $K=\overline{f(\R^d)}$ we obtain the global Wiener--Levi Theorem for functions from the class $W$.

 Note that exponents of $g$ belong to the span over $\Z$ of exponents of $f$.
\medskip

 The main consequence  of Theorem \ref{T} is the following result:
 \begin{Th}[\cite{F0}]\label{1}
For every $f\in W$ and $\e>0$ there is a function $g\in W$ such that $f(x)g(x)=1$ for all $x\in\R^d$ such that $|f(x)|\ge\e$ and $g(x)=0$ for all $x\in\R^d$ such that $|f(x)|\le\e/2$.
  \end{Th}

{\bf 2. Tempered distributions and measures and their Fourier transform.}
 To show applications of Theorem \ref{1}, we recall some definitions (see, for example, \cite{R}).

Denote by $S(\R^d)$ the Schwartz space of test functions $\p\in C^\infty(\R^d)$ with finite norms
 $$
  p_m(\p)=\sup_{\R^d}(\max\{1,|x|\})^m\max_{k_1+\dots+k_d\le m} |D^k(\p(x))|,\quad m=0,1,2,\dots,
 $$
  $k=(k_1,\dots,k_d)\in(\N\cup\{0\})^d,\  D^k=\partial^{k_1}_{x_1}\dots\partial^{k_d}_{x_d}$. These norms generate a topology on $S(\R^d)$, and elements of the space $S^*(\R^d)$ of continuous linear functionals on $S(\R^d)$ are called tempered distributions.
  For every tempered distribution $f$ there exist $C>0$ and $m\in\N\cup\{0\}$ such that for all $\p\in S(\R^d)$
$$
                           |f(\p)|\le Cp_m(\p).
$$

The Fourier transform of a tempered distribution $f$ is defined by the equality
\begin{equation}\label{f}
\hat f(\p)=f(\hat\p)\quad\mbox{for all}\quad\p\in S(\R^d),
\end{equation}
where
$$
   \hat\p(y)=\int_{\R^d}\p(x) e^{-2\pi i\la x,y\ra}dx
 $$
is the Fourier transform of the function $\p$. Note that the Fourier transform of every tempered distribution is also a tempered distribution. But here we consider only the case when $f$ and $\hat f$ are complex Radon measures on $\R^d$. For example, if
$$
f(x)=\sum_n c_ne^{2\pi i\la x,\g_n\ra}\in W,
$$
then the Fourier transform of the measure $f(x)dx$ is equal to $\sum_n c_n\d_{\g_n}$, where $\d_\g$ means the unit mass at $\g\in\R^d$.
Also, if $\mu^0=\sum_{k\in\Z^d}\d_k$, then by the Poisson formula,
$$
  \sum_{n\in\Z^p}f(n)=\sum_{n\in\Z^p}\hat f(n), \quad f\in S(\R^d),
$$
 and we have $\hat\mu^0=\mu^0$.
Therefore,  if $L$ is a full-rank lattice, i.e., $L=A(\Z^d)$ for some nonsingular linear operator $A$ in $\R^d$, and $\mu^1=\sum_{\l\in L+a}\d_\l$ for some $a\in\R^d$, then
\begin{equation}\label{p}
\hat\mu^1(dy)=(\det A)^{-1}\sum_{\l\in L^*} e^{-2\pi i\la y,a\ra}\d_\l(dy),
\end{equation}
where $L^*=\{y\in\R^d: <\l,y>\in\Z\quad \forall \l\in L\}$ is the conjugate  lattice.
\medskip

 Let $\nu$ be a Radon measure on $\R^d$. Denote by $\nu_t$ its shift on $t\in\R^d$, i.e.,
 $$
 \int g(x)\nu_t(dx)=\int g(x+t)\nu(dx).
 $$
  A  measure $\nu$ is  {\it translation bounded} if variations of its translations $|\nu_t|$  are  bounded in the unit ball  uniformly in $t\in\R^d$. Note that every translation bounded measure on $\R^d$ satisfies the condition
 \begin{equation}\label{m}
 |\nu|(B(0,r))=O(r^d),\quad  r\to\infty,
 \end{equation}
 therefore it belongs to $S^*(\R^d)$. Here $B(x,r)=\{t\in\R^d:\,|t-x|<r\}$.

 The measure $\nu$ on $\R^d$ is {\it atomic}, if it has the form
 $$
 \nu=\sum_{\l\in\L} c_\l\d_\l,\quad c_\l\in\C,\ \mbox{ with countable }\L\subset\R^d.
$$
 If this is the case, we will write $\nu(\l):=c_\l$. Also, we shall say that $\L$ is a {\it support} of $\nu$.

 The measure $\nu$ on $\R^d$ has {\it a uniformly discrete support $\L$}, if
 $$
 \inf\{|x-x'|:\,x,\,x'\in\L, x\neq x'\}>0.
 $$

Such measures are the main object in the theory of Fourier quasicrystals (see \cite{C}-\cite{F2}, \cite{F0}-\cite{La1}, \cite{LO1}-\cite{Mo1}). This theory was developed in connection with the  experimental discovery  of non-periodic atomic structures with diffraction patterns consisting of spots, which was made in the mid '80s.

Remark also that some properties of tempered distributions with discrete and closed or atomic supports were considered in \cite{F3} - \cite{F0}, \cite{LA}, \cite{P}. In particular, the following statements are implicitly contained in \cite{F0}. But for completeness we present proofs  of them at the end of this article.
\medskip

\begin{Pro}\label{P3}
If $\nu$ is a translation bounded measure and $\psi\in S(\R^d)$, then the total mass of the variation $|\psi\nu_t|$ of the measure $\psi\nu_t$ is  bounded uniformly in $t\in\R^d$. Moreover, for every $\e>0$ there is $r(\e)<\infty$ such that the mass of restriction of each measure $|\psi\nu_t|$ on the set $\{x\in\R^d:\,|x|>r(\e)\}$ is less then $\e$.
\end{Pro}
\begin{Pro}[also, see \cite{KL} and \cite{LA}]\label{P0}
     If $\nu$ is a  measure from $S^*(\R^d)$ with uniformly discrete support $\L$, and $\hat\nu$ is a measure satisfying (\ref{m}), then $\sup_{\l\in\L}|\nu(\l)|<\infty$, hence the measure $\nu$ is translation bounded.
\end{Pro}
\begin{Pro}\label{P1}
If $\nu$ is a measure from $S^*(\R^d)$, $\hat\nu$ is an atomic measure satisfying (\ref{m}), and $\psi\in S(\R^d)$,  then the convolution $(\psi\star\nu)(t)=\int\psi(t-x)\nu(dx)$ belongs to $W$, and its Fourier transform equals $\hat\psi\hat\nu$.
\end{Pro}
\begin{Pro}\label{P2}
If $\nu$ is a translation bounded measure,  $\hat\nu$ is a translation bounded atomic measure, and $g\in W$, then the Fourier transform  $\widehat{g\nu}$ of the product $g\nu$ is a translation bounded atomic measure.
  \end{Pro}
\medskip

 {\bf 3. Properties of measures with uniformly discrete support.} The following theorem belongs to Y.Meyer \cite{M1}.
\begin{Th}
Let $\mu=\sum_{\l\in\L} a_\l\d_\l$ be a  measure on $\R$ with a discrete and closed support $\L$ and the set $\{x=a_\l:\,\l\in\L\}$ is finite.
If $\mu\in S^*(\R)$ and its Fourier transform $\hat\mu$ is a translation bounded measure on $\R$, then
 $$
 \L=E\triangle\bigcup_{j=1}^N(\a_j\Z+\b_j),\quad \a_j>0,\ \b_j\in\R,
 $$
where the set $E$ is finite.
\end{Th}
Here $A\triangle B$ means the symmetrical difference between $A$ and $B$.
\medskip

 In \cite{K} M.Kolountzakis extended the above theorem to measures on $\R^d$. Also, he replaced the condition "the measure $\hat\mu$ is translation bounded" with the weaker one "the measure $\hat\mu$ satisfies condition (\ref{m})". He also found a condition for the support of $\mu$ to be a finite union of translations of several full-rank lattices. His result is very close to Cordoba's  one:
  \begin{Th}[\cite{C}]
 Let  $\mu=\sum_{\l\in\L} a_\l\d_\l$ be a measure on $\R^d$ with a uniformly discrete support $\L$ and a finite set $\{x=a_\l:\,\l\in\L\}$.  If $\hat\mu$ is atomic and translation bounded measure,  then $\L$  is a finite union of translations of several, possibly incommensurable, full-rank lattices.
\end{Th}
\noindent In papers \cite{F2} and \cite{F3} we get a small elaboration of Cordoba's result. In particular, we replaced the conditions "$a_\l$ from a finite set" by "$|a_\l|$ from a finite set". But Cordoba's type theorems are not true for some uniformly discrete measures $\mu=\sum_{\l\in\L} a_\l\d_\l$ with translation bounded $\hat\mu$ and a countable set $\{a_\l\}_{\l\in\L}$ (\cite{LO2}).
\medskip

 Theorem \ref{1}  makes it possible to obtain such a result:
\begin{Th}[\cite{F0}]\label{2}
 Let $\mu=\sum_{\l\in\L} a_\l\d_\l$  be a measure on $\R^d$ with  a uniformly discrete support $\L$ such that $\inf_\L|a_\l|>0$,  and $\hat\mu$  be an atomic measure satisfying (\ref{m}). Then  $\L$  is a finite union of translations of several disjoint full-rank lattices.
 \end{Th}
Note that even if both $\supp\mu$ and $\supp\hat\mu$ are uniformly discrete, they can be  finite unions of translations of incommensurate lattices (\cite{F2}).

\medskip
We supplement Theorem \ref{2} with a description of the measure $\mu$.

\begin{Th}\label{3}
 In conditions of Theorem \ref{2},
\begin{equation}\label{r1}
\mu=\sum_{j=1}^N F_j(y)\D^j,
\end{equation}
 where $\D^j$ are sums of unit masses at the points of some lattices $L_j$ or their translations,
and $F_j(y)=\sum_s b_s e^{2\pi i\la y,\a_s^j\ra}\in W$ with a bounded set of $\a_s^j\in\R^d$.
Moreover,
 \begin{equation}\label{r2}
 \hat\mu=\sum_{j=1}^N e^{2\pi i\la x,\l_j\ra}\nu^j,
 \end{equation}
  where $\nu^j$ are $d$-periodic atomic measures with full-rank lattices $L_j^*$ of periods,
 and $\l_j\in\L$.
\end{Th}

{\bf Proof of Theorem \ref{3}}. Let $\eta<(1/2)\inf\{|\l-\l'|:\, \l,\l'\in\G\}$, and $\psi$ be an odd $C^\infty$ function such that $\psi(0)=1$ and $\supp\psi\subset B(0,\eta)$. Clearly, $(\psi\star\mu)(\l)=\mu(\l)$ for $\l\in\L$\ and,
by Proposition \ref{P1}, the function $g=\psi\star\mu\in W$.

By Theorem \ref{2}, the support $\L$ of the measure $\mu$ is a finite union of translations of  full-rank lattices $\l_j+L_j, \ j=1,\dots,N$. Set
$$
\D^j=\sum_{\l\in L_j+\l_j}\d_\l.
$$
We have
\begin{equation}\label{h}
 \mu=\sum_{j=1}^N g\D^j,
 \end{equation}
  with $g(x)=\sum_n c_ne^{2\pi i\la x,\g_n\ra}\in W$. For every fixed $j$ and each $\g\in\R^d$ there is $\a$ inside the parallelepiped generated by corresponding $L_j^*$ such that $\g-\a\in L_j^*$, therefore,
$e^{2\pi i\la x,\g\ra}=e^{2\pi i\la x,\a\ra}$ for $x\in L_j$.  Collecting similar terms for this $j$, we obtain (\ref{r1}).

 Next, by (\ref{p}),  $\hat\D^j$ is a uniformly discrete and translation bounded measure. Furthermore, the measure $g\D^j$ coincides with the restriction of the measure $\mu$ to $\l_j+L_j$ and, by Proposition \ref{P2}, its  Fourier transform $\widehat{g\D^j}$ is an atomic translation bounded measure. Set
$$
\nu^j=e^{-2\pi i\la\l_j,y\ra}\widehat{g\D^j}.
$$
The converse Fourier transform of $\nu^j$ is equal to
$$
(\nu^j)\check\ =\sum_{x\in L_j+\l_j}\mu(x)\d_{x-\l_j}=\sum_{x\in L_j}\mu(x+\l_j)\d_x,
$$
 and the converse Fourier transform of the measure $\nu^j_a$ for $a\in\R^d$ is equal to
 $$
 e^{2\pi i\la a,x\ra}\sum_{x\in L_j}\mu(x+\l_j)\d_x,
 $$
  which coincides with $(\nu^j)\check\ $ for each $a\in L_j^*$.
 Therefore, $\nu^j_a=\nu^j$. So, $\nu^j$ is  $d$-periodic with the  lattice $L_j^*$ of periods and, by (\ref{h}),
 $$
 \hat\mu=\sum_{j=1}^N\widehat{g\D^j}=\sum_{j=1}^N e^{2\pi i\la\l_j,y\ra}\nu^j.
 $$
 \bs
\medskip

{\bf 4. Coherent sets of frequencies.}
Let us remember that a uniformly discrete set $\U\subset\R^d$ is  {\it a coherent set of frequencies} (or satisfies Kahane's property), if every limit of a sequence of finite sums
$$
\sum c_\l e^{2\pi i\la x,\l \ra},\quad\l\in\U,\quad c_\l\in\C,
$$
with respect to the topology of uniform convergence on every compact subset of $\R^d$ is almost periodic in the sense of H.Bohr.
 \begin{Th}[Y.Meyer, \cite{M2}]\label{3a}
 Let $\mu=\sum_{\l\in\L} a_\l\d_\l$ be a Radon measure from $S^*(\R^d)$ with a uniformly discrete support $\L$ and $\hat\mu$ be a translation bounded Radon measure.
  Then the set $\U=\{\l:\, a_\l=1\}$ is a coherent set of frequencies.
  \end{Th}

Using Theotem \ref{1}, we obtain the following result.
\begin{Th}\label{4}
 Let measures $\mu=\sum_{\l\in\L} \mu(\l)\d_\l$ and $\hat\mu$ be atomic translation bounded measures from $S^*(\R^d)$, and let $\U\subset\L$ be such that for all $\l\in\U$ and some $\e>0$

   i) $|\mu(\l)|\ge\e$,

   ii) $|\l-\l'|\ge\e$ for all $\l'\in\L\setminus\{\l\}$.

\noindent Then the set $\U$ is a coherent set of frequencies.
\end{Th}
  \begin{Cor}
 Let $\mu=\sum_{\l\in\L} \mu(\l)\d_\l$ be a  measure from $S^*(\R^d)$ with a uniformly discrete support $\L$ and $\hat\mu$ be an atomic translation bounded  measure.
     Then for every $\e>0$ the set $\U=\{\l:\, |\mu(\l)|\ge\e\}$ is a coherent set of frequencies.
  \end{Cor}
  Indeed, by Proposition \ref{P0},  the measure $\mu$ is translation bounded, therefore  all the conditions of Theorem \ref{4} are met.
  \medskip

{\bf Proof of Theorem \ref{4}}.
Let $\eta<\e/2$, and $\psi$ be the same as in the proof of Theorem \ref{3}. By Proposition \ref{P1}, the function $g=\psi\star\mu\in W$.
Then $g(\l)=(\psi\star\mu)(\l)=\mu(\l)$ for $\l\in\U$.  By Theorem \ref{1}, there is $h\in W$ such that $g(\l)h(\l)=1$ under condition $|g(\l)|\ge\e$, in particular, for all $\l\in\U$. Fix a parameter $t\in\R^d$. Let
$F(x)$ be a convolution of the function $\psi$ and the measure $e^{2\pi i\la x, t\ra}h(x)\mu(x)$, i.e.,
$$
F(x)=\sum_{\l\in\L}\psi(x-\l)e^{2\pi i\la\l, t\ra}h(\l)\mu(\l).
$$
By  Proposition \ref{P2}, the measure $\widehat{h\mu}$ is translation bounded. Using Proposition \ref{P1}, we see that
$\hat F(y)=\hat\psi(\widehat{h\mu})_t(y)$. Applying Proposition \ref{P3}, we get that the total mass of the measure $\hat F(y)$  is bounded by some constant $C<\infty$, and the mass of its restriction  to the set $\{x\in\R^d:\,|x|>r\}$  is less than $1/2$ for a suitable  $r<\infty$. Note that $C$ and $r$ are independent of $t$. Taking into account that $F(x)$ is the converse Fourier transform of the measure $\hat F(y)$ and obvious equality
$F(\l)=e^{2\pi i\la\l, t\ra}$ for all $\l\in\U$, we get
$$
  e^{2\pi i\la\l, t\ra}=\int_{\R^d} e^{2\pi i\la y,\l\ra}\hat F(dy).
 $$
Now, let $\sum c_\l e^{2\pi i\la t,\l\ra}$ be any finite sum of exponents with $\l\in\U$. We have
$$
   \left|\sum c_\l e^{2\pi i\la t,\l\ra}\right|= \left|\int_{\R^d}\sum c_\l e^{2\pi i\la y,\l\ra}\hat F(dy)\right|
$$
$$
 \le \left|\int_{|y|\le r}\sum c_\l e^{2\pi i\la y,\l\ra}\hat F(dy)\right|+\left|\int_{|y|>r}\sum c_\l e^{2\pi i\la y,\l\ra}\hat F(dy)\right|.
 $$
 The first integral in the right-hand side does not exceed
 $$
 C\sup_{|y|\le r}\left|\sum c_\l e^{2\pi i\la y,\l\ra}\right|,
 $$
 the second one is bounded by
 $$
 \frac{1}{2}\sup_{t\in\R^d}\left|\sum c_\l e^{2\pi i\la t,\l\ra}\right|.
 $$
 Thus,
 $$
 \sup_{t\in\R^d}\left|\sum c_\l e^{2\pi i\la t,\l\ra}\right| \le 2C\sup_{|y|\le r}\left|\sum c_\l e^{2\pi i\la y,\l\ra}\right|.
 $$
 Therefore, the uniform convergence of the sequence of exponential sums
on the ball $B(0,r)$ implies the uniform convergence on $\R^d$, and the limit of the sequence is an almost periodic function in the sense of H.Bohr. \bs
\medskip

{\bf 5. Proofs of propositions \ref{P3} - \ref{P2}}. To prove Proposition \ref{P3}, fix $t\in\R^d$. Denote by $N(r)$ the variation $|\nu_t|$ in the ball $B(0,r)$. Since the measure $\nu$ is translation bounded, we see that $N(r)\le C(1+r)^d$ with a constant $C$ independent of $t$.  Also, $|\psi(x)|\le C'(1+|x|)^{-d-1}$ for all $x\in\R^d$ with a constants $C'$. Therefore, integrating by parts, we obtain the estimate
 $$
   \int_{\R^d}|\psi(x)||\nu_t|(dx)\le C'\int_0^\infty(1+r)^{-d-1}dN(r)\le CC'(d+1)\int_0^\infty(1+r)^{-2}dr.
 $$
 Also, if $r(\e)$ is sufficiently large, we get
 $$
   \int_{|x|>r(\e)}|\psi(x)||\nu_t|(dx)\le C'\int_{r(\e)}^\infty(1+r)^{-d-1}dN(r)\le CC'(d+1)\int_{r(\e)}^\infty(1+r)^{-2}dr<\e.
 $$
 \medskip

To prove  Proposition \ref{P0}, note that (\ref{f}) implies the equality
 $$
   \nu(\l)=\int\hat\psi(x-\l)\nu(dx)=\int\psi(y)e^{2\pi i\la\l,y\ra}\hat\nu(dx),
 $$
 where $\psi$ is the same function as in the proof of Theorem \ref{3}. Arguing as in the proof of Proposition \ref{P3}, we see that the module of the latter integral is bounded uniformly in $\l\in\L$.
  \medskip

  To prove  Proposition \ref{P1}, we can repeat  arguments from the proof of Proposition \ref{P3} and get that the total variation of the measure $\hat\psi\hat\nu$ is finite. Since the measure  $\hat\psi\hat\nu$ is atomic, we
   see that it has  the form
  $$
  \sum_{n=1}^\infty b_n\d_{\g_n},\qquad \sum_n|b_n|<\infty,
  $$
 and its converse Fourier transform is equal
 $$
  \sum_{n=1}^\infty b_n e^{2\pi i\la \g_n,x\ra}\in W.
  $$
On the other hand, the converse Fourier transform of the function $\p(x):=\psi(x-t)\in S(\R^d)$ equals $\hat\psi(y)e^{2\pi i\la y,t\ra}$. Since $\psi$ is an odd function, we get by formula (\ref{f})
 $$
   (\psi\star\nu)(t)=\int\psi(t-x)\nu(dx)=\int e^{2\pi i\la t,y\ra}\hat\psi(y)\hat\nu(dy).
 $$
Therefore, the converse Fourier transform of the measure $\hat\psi\hat\nu$ is equal $\psi\star\nu$.
\medskip

To prove Proposition \ref{P2}, note that the function $g$ is bounded, hence the measure $g\nu$ is translation bounded. Set
 $$
 \nu^\g(x)=e^{2\pi i\la x,\g\ra}\nu(x),\quad \g\in\R^d.
 $$
 Suppose $g(x)=\sum_n c_ne^{2\pi i\la x,\g_n\ra}\in W$. Then
 $$
 g\nu=\sum\nolimits_n c_n\nu^{\g_n},\qquad  \widehat{g\nu} =\sum\nolimits_n c_n\widehat{\nu^{\g_n}}.
 $$
 Note that $\sum_n|c_n|<\infty$, and $\widehat{\nu^{\g_n}}$ are atomic measures, hence $\widehat{g\nu}$  is an atomic measure too.
  Then for each $y\in\R^d$  we get $|\widehat{\nu^\g}|(B(y,1))=|\hat\nu|(B(y-\g,1))$. Therefore,
$$
|\widehat{g\nu}|(B(y,1))\le \sum_n |c_n||\widehat{\nu^{\g_n}}|(B(y,1))\le \sup_{t\in\R^d}|\hat\nu|(B(t,1))\sum_n |c_n|.
$$
 So, $\widehat{g\nu}$ is a translation bounded atomic measure. \bs
 
\end{document}